\documentclass[11pt]{article} \usepackage{amssymb}
\usepackage{amsfonts} \usepackage{amsmath} \usepackage{amsthm} \usepackage{bm}
\usepackage{latexsym} \usepackage{epsfig}
\usepackage{hyperref} 

\newtheorem*{theorem*}{Theorem}
\newtheorem{theorem}{Theorem}[section]
\newtheorem{lemma}[theorem]{Lemma}
\newtheorem*{proposition*}{Proposition}

\newtheorem{claim}[theorem]{Claim}

\newtheorem{definition}[theorem]{Definition}

\newcommand{\ignore}[1]{}

\newcommand{\enote}[1]{} \newcommand{\knote}[1]{}
\newcommand{\rnote}[1]{}

\renewcommand{\P}[1]{{\mathbb{P}}\left[{#1}\right]}
\newcommand{\CondP}[2]{{\mathbb{P}}\left[{#1}\middle\vert{#2}\right]}

\newcommand{\ind}[1]{{\bf 1}\left({#1}\right)}

\newcommand{\CondE}[2]{{\mathbb{E}}\left[{#1}\middle\vert{#2}\right]}

\newcommand{\N}{\mathbb N} \newcommand{\R}{\mathbb R}

\newcommand{\half}{{\textstyle \frac12}}

\newcommand{\argmax}{\operatornamewithlimits{argmax}}

\renewcommand{\phi}{\varphi}

\newcommand{\vote}{V}
\newcommand{\votevec}{\bar{\vote}}
\newcommand{\psignal}{W}

\begin{document}


\title{Making Consensus Tractable}

\author{Elchanan Mossel\footnote{UC Berkeley and Weizmann Institute of
    Science. Supported by a Sloan fellowship in Mathematics, by BSF
    grant 2004105, by NSF Career Award (DMS 054829) by ONR award
    N00014-07-1-0506 and by ISF grant 1300/08.} and Omer
  Tamuz\footnote{Weizmann Institute of Science. Omer Tamuz is
    supported by ISF grant 1300/08, and is a recipient of the Google
    Europe Fellowship in Social Computing. This research is supported
    in part by this Google Fellowship.} }

\maketitle

\begin{abstract}
  We study a model of consensus decision making, in which a finite
  group of Bayesian agents has to choose between one of two courses of
  action. Each member of the group has a private and independent
  signal at his or her disposal, giving some indication as to which
  action is optimal. To come to a common decision, the participants
  perform repeated rounds of voting. In each round, each agent casts a
  vote in favor of one of the two courses of action, reflecting his or
  her current belief, and observes the votes of the rest.
  
  We provide an efficient algorithm for the calculation the agents
  have to perform, and show that consensus is always reached and that
  the probability of reaching a wrong decision decays exponentially
  with the number of agents.

\end{abstract}

\section{Introduction}
Consensus voting, or decision by unanimous agreement, is a method of
communal governance that requires all members of a group to agree on a
chosen course of action. The European Union's Treaty of
Lisbon~\cite{LisbonTreaty:07} decrees that {\em ``Except where the
  Treaties provide otherwise, decisions of the European Council shall
  be taken by consensus.''} In this the EU follows the historical
example of the Diet of the Hanseatic League~\cite{Miller:09} and
others.

Proponents of this method consider it to have many advantages over
majority voting: it cultivates discussion, participation and
responsibility, and avoids the so-called ``tyranny of the
majority''. The drawback is, of course, a lengthy and difficult decision
making process, lacking even the guarantee of a conclusive ending.

However, in standard theoretical setups of rational Bayesian
participants (e.g. \cite{McKelvey:86}, \cite{GaleKariv:03}),
agents cannot ``agree to disagree''~\cite{Aumann:76}, and
consensus is eventually reached. Unfortunately, this may come at the
price of tractability; Bayesian calculations can, in some situations,
be practically impossible~\cite{GaleKariv:03}.

Indeed, modeling economic behavior involves an inherent conflict
between rationality and tractability~\cite{Simon:82}. It seems that in
many situations it is practically impossible to calculate which course
of action is optimal, leaving the theoretician with a model that is
either not rational, and thus hard to justify, or not tractable and
hence not realistic. A common course of action is to relax the
rationality assumption and consider boundedly-rational agents. We do
not do this, but rather show that a fully rational model is, perhaps
unexpectedly, tractable.

We consider a model describing a group of Bayesian agents that have to
make a binary decision. We show that under the dynamics we describe,
unanimity is reached with probability one, and give an efficient
algorithm for the agents' calculations.

Our model features a finite group of Bayesian agents that have to
choose between two possible courses of action. Each initially receives
a private and independent signal, which contains some information
indicating which action is more likely to be the correct one. The
agents participate in rounds of voting, in which each indicates which
action it believes is more likely to be correct, and learns the
others' opinion thereof. The process continues until unanimity is
reached. The Bayesian agents are myopic, so their actions are not
strategic, but truthfully reflect the information available to
them. They are rational and do not follow heuristics, or rules of
thumb, or boundedly-rational courses of action.

As an example, consider a committee that has to decide whether or not
to accept a candidate for a position, who a-priori has a chance of one
half to be a good hire. Each committee member gets to interview the
candidate in private. If the candidate is good - i.e., the correct
action is ``hire'', then each committee member $i$ receives a private
signal $\psignal_i$ drawn independently from $N(1,1)$, the normal
distribution with expectation 1 and variance 1. If the candidate is
bad (i.e., the correct action is ``don't hire''), then $\psignal_i$ is
drawn from $N(-1,1)$.

The committee now commences to vote in rounds. In each round of voting
each member casts a public ``hire'' or ``don't hire'' vote, depending
on which it thinks is more likely to be the correct decision. At each
iteration, each member's opinion is based on its private signal, as
well as the votes of the other members in the previous rounds.

The agents are Bayesian in the sense that their beliefs are precisely
calculated according to Bayes' Law. This is not a straightforward
calculation, as they have to take into account that each of the votes
of their peers was also likewise calculated.  From the description of
the process it is not at all clear that there is a succinct
description of the decision process taken by the agents, and how can
this process be analyzed (indeed - we challenge the reader to try!).

In slightly more general settings the problem seems even more
difficult: for example, consider the case that the agents lie on a
social network graph, i.e., when each agent only observes the actions
of only a subset of the rest. There, no efficient algorithms for the
agents' calculations are known, and it is in fact conjectured that
none exist (cf., Kanoria and Tamuz~\cite{KanoriaTamuz12}).  When more
than two possible actions are available, then too it seems that the
computational problem is significantly more difficult, although
perhaps not intractable. An interesting open problem is to suggest an
efficient algorithm for the agents' calculation in these (more
general) models.

For our model we show that there does exist an efficient algorithm to
perform the agents' calculations. We also show that the agents will
eventually all cast the same vote, and that the probability that this
vote is correct approaches one as the number of agents increases.

\subsection{Related work}
In 1785 the Marquis de Condorect proved a founding
result~\cite{Condorcet:85} in the field of group decision making. The
Condorect Jury Theorem states that given that each member of a jury
knows the correct verdict with some probability $p>1/2$, the
probability that the jury reaches a correct decision by a majority
vote goes to ones as the size of the jury increases. Our ``asymptotic
learning'' result extends this theorem to a more general class of
private signals, given that at least two rounds of voting are carried
out.

Sebenius and Geanakoplos~\cite{sebenius1983don} show in a classical
paper that a pair of agents eventually reach agreement on the ``state
of with world'' in a model similar to ours, with finite probability
spaces. Likewise, a consequence of the convergence proof given by Gale
and Kariv~\cite{GaleKariv:03}\footnote{See a comment on this paper in
  ~\cite{mueller2010comment}.} as well as Rosenberg, Solan and
Vieille~\cite{rosenberg2009informational} and
Mueller-Frank~\cite{mueller2011general} (three models which are
generalizations of ours) is that if a pair of agents' actions
converge, it is to the same action, unless the agents are indifferent
at the limit. However, none of these results imply that the agents'
actions do in fact converge, or that the agents reach agreement.

We provide a basic proof of a stronger result, namely that all the
agents' actions converge, and in particular to the same action. We
further show that each round of voting increases the probability that
a given agent votes for the better alternative, and that this
probability goes to one at the second round of voting, as the number
of participants goes to infinity.  Finally, our most significant
improvement over the work of Gale and Kariv is that we provide the
participants with an efficient algorithm to calculate their beliefs.

The question of complexity of agreement was discussed from a somewhat
different perspective in an interesting paper by
Aaronson~\cite{aaronson2005complexity}. Aaronson discusses the
complexity of agreement in a revealed beliefs model (where agents
reveal their beliefs rather than just observe actions), but with
general correlated signals. He provides polynomial bounds for
approximate convergence (in the sense that the beliefs are close at a
specific time $t$, but with no guarantee on how close they are at
later times) and designs even more efficient algorithms for achieving
agreement. Of course, in the context of the current paper, if the
agents reveal their beliefs, then consensus is achieved in one
round. While our paper is much more restricted in terms of the graph
(the complete graph) and the signal structure (conditional
independence) it considers a natural action dynamics, as opposed to
Aaronson's belief dynamics or artificial dynamics. Our results are
also stated in terms of convergence of beliefs and not just closeness
of beliefs at a certain point in time.

In a subsequent work to this paper, together with Allan
Sly~\cite{mossel2012asymptotic}, we show that for very general
voting models asymptotic learning holds in the sense that as the
number of voters goes to infinity the probability of convergence to
the correct outcome goes to one. However, the results of this article
are not known to extend to the general models studied
in~\cite{mossel2012asymptotic}. There, no efficient update
algorithms are known, no rates of convergence are known, and it is not
known whether agents always reach consensus. We answer all of these
question for the model we study.

Note that the results of~\cite{mossel2012asymptotic} do imply that
as the number of agents goes to infinity, the probability of
non-convergence (or convergence to the wrong state) goes to
zero. However, as in the work of Gale and Kariv, for finite graphs it
remains possible that with positive probability both actions are taken
infinitely often, with beliefs converging to values implying
indifference between the actions. We show that this is not the case in
our model.

Studies in committee mechanism design
(cf.~\cite{laslier2008committee,glazer1998motives}) strive to
construct mechanisms for eliciting information out of committee
members to arrive at optimal results. We, by and large, do not take
this path but consider a ``natural'' setting which was not
specifically constructed to achieve any such goal. 

One could indeed raise the objection that the process could be made
simpler if the committee members were to tell each other their private
signals, in which case the optimal answer would be arrived at
immediately. However, a common assumption in the study of Bayesian
agents (cf.~\cite{Banerjee:92,BichHirshWelch:92,smith2000pathological,GaleKariv:03}) is that ``actions speak louder
than words'', so that agents learn from each other's actions rather
than revealing to each other all their information. The latter option
may not be feasible, as the said information may consist of
experiences and impressions that could take too long to explain, may
be difficult to articulate, or may not be even consciously known. In
our case the agents' actions are the casting of votes.

Our work is more closely related to models of herd behavior
(cf.~\cite{Banerjee:92,BichHirshWelch:92,smith2000pathological}). These
feature a group of agents with a state of the world and corresponding
private signals, much like ours. There too agents don't observe each
others' private signals but only actions. However, there the agents
are exogenously ordered, and each takes a single action after seeing -
and learning from - the actions of its predecessors.

\subsection{Model}
Our model features a finite set of agents $[n]=\{1,2,\ldots,n\}$ that have
to make a binary decision regarding an unknown state of the world $S
\in \{0,1\}$. Each is initially given a private signal $\psignal_i$,
distributed $\mu_0$ if $S=0$ and $\mu_1$ if $S=1$, and independent of the
other signals, conditioned on $S$.

\begin{definition}
  \label{def:assumptions}
  Let $\mu_0$ and $\mu_1$ be measures on a $\sigma$-algebra
  $(\Omega,\mathcal{O})$ satisfying the following conditions:
  \begin{enumerate}
  \item $\mu_0$ and $\mu_1$ are mutually absolutely continuous, so
    that, by the Radon-Nikodym theorem, the Radon-Nikodym derivative
    $\frac{d\mu_1}{d\mu_0}(\omega)$ exists and is non-zero for all
    $\omega\in\Omega$.
  \item \label{item:non-atomic} Let $\psignal$ be distributed $\half \mu_0 + \half \mu_1$, and
    let $X = \log\frac{d\mu_1}{d\mu_0}(\psignal)$. Then the
    distribution of $X$ is non-atomic. 
  \end{enumerate}
\end{definition}
Note that (\ref{item:non-atomic}) implies that $\mu_0 \neq \mu_1$,
since otherwise $X=0$ a.s.\ and thus its distribution is atomic.

\begin{definition}
  \label{def:signals}
  Let $\mu_0$ and $\mu_1$ be measures on a $\sigma$-algebra
  $(\Omega,\mathcal{O})$ satisfying the conditions of
  definition~\ref{def:assumptions}. Let $\delta_0,\delta_1$ denote the
  measures on $\{0,1\}$ that satisfy $\delta_0(0)=\delta_1(1)=1$ and
  $\delta_0(1)=\delta_1(0)=0$.

  Let $\mathbb{P}$ be the probability measure over the space
  $\{0,1\}\times\Omega^n$ given by
  \begin{align}
    \label{eq:P-def}
    \mathbb{P} = \half\delta_0\otimes{\mu_0}^{\otimes n} +
    \half\delta_1\otimes{\mu_1}^{\otimes n}.
  \end{align}
  Let $S$, taking values in $\{0,1\}$, and
  $(\psignal_1,\ldots,\psignal_n)$, taking values in $\Omega^n$, be
  random variables with joint distribution $\mathbb{P}$:
  \begin{align*}
  (S, \psignal_1, \ldots, \psignal_n) \sim \mathbb{P}.    
  \end{align*}
  We call $S$ the {\em state of the world} and call $\psignal_i$ agent
  $i$'s {\em private signal}.
\end{definition}

%

Equivalently, $S$ is picked uniformly from $\{0,1\}$, and conditioned
on $S$, the agents' private signals $\psignal_i$ are distributed
i.i.d.\ $\mu_S$. It is important to note that conditioned on $S$ the
private signals $\psignal_i$ are independent. In much of what follows
it is not necessary to assume that they are identical, but we make
this assumption to simplify notation and conform to a widely studied
economic model.

The agents participate in a process of voting rounds. In each round
$t$ each agent $i$ casts a public vote $V_i(t) \in \{0,1\}$, depending
which of the two is more likely to be the state of the world,
conditioned on the information available to $i$; this includes
$\psignal_i$ as well as the votes of the other agents in the previous
rounds.

\begin{definition}
  \label{def:voting}
  For $t \in \{1,2,\ldots\}$ and agent $i \in [n]$, let $\vote_i(t)$, the
  {\em vote} of agent $i$ at time $t$, be defined by
  \begin{align}
    \label{eq:V-def}
    \vote_i(t) =
    \begin{cases}
      1 & \mbox{if } \CondP{S=1}{\psignal_i,\votevec^{t-1}} > 1/2 \\
      0 & \mbox{otherwise} \\
    \end{cases}
  \end{align}
  where $\votevec^t = \{\vote_j(t'):\: j \in [n],\, t' \leq t\}$ denotes the
  votes of all agents up to time $t$.
\end{definition}
Alternatively, one could define
\begin{align}
  \label{eq:vote-max-prob}
  \vote_i(t) = \argmax_{s \in \{0,1\}}\CondP{S=s}{\psignal_i,\votevec^{t-1}},
\end{align}
with a ``tie breaking law'' specifying that when the conditional
probability on the r.h.s.\ is half then the vote is 0.  We note that
it is easy to see that the assumption that the distribution of the
Radon-Nikodym derivatives $\frac{d\mu_1}{d\mu_0}(\psignal_i)$ is
non-atomic (definition~\ref{def:assumptions}) implies that the
probability of ever encountering a tie is $0$ a.s.\ and therefore the
details of the tie breaking rule a.s.\ do not affect the behavior of
the process or our results.

\subsection{Results}

For the model defined in definitions~\ref{def:assumptions},
\ref{def:signals} and~\ref{def:voting} we prove the
following theorems.
\begin{itemize}
\item {\bf Unanimity:} A unanimous decision is always reached. That
  is, with probability one all agents vote identically at some
  round, and the process essentially ends.
\begin{theorem}\label{thm:unanimity}
  With probability one there exists a time $T_u$ and a vote $\vote \in
  \{0,1\}$ such that for all $t \geq T_u$ and agents $i$ it holds that
  $\vote_i(t)=\vote$.
\end{theorem}
\item {\bf Monotonicity:} The probability that an agent votes
  correctly is non-decreasing with the progression of rounds.
\begin{theorem}
  \label{thm:monotonicity}
  For all agents $i$ and times $t>1$, it holds that
  \begin{align*}
    \P{\vote_i(t)=S} \geq \P{\vote_i(t-1)=S}.
  \end{align*}
\end{theorem}
\item {\bf Asymptotic Learning:} The probability of reaching a correct
  decision at the end of the process approaches one as the number of
  agents increases. In fact, this already holds by the second
  round of voting.
\begin{theorem}
  \label{thm:learning}
  Fix $\mu_0$ and $\mu_1$, and let $n$ be the number of agents.  Then
  there exist constants $C=C(\mu_0,\mu_1)$ and $n_0=n_0(\mu_0,\mu_1)$ such
  that
  \begin{align*}
    \P{\forall i: \vote_i(2) = S} > 1-e^{-Cn}
  \end{align*}
  for all $n > n_0$.
\end{theorem}
\item {\bf Tractability:} In order to discuss tractability we must
  assume that certain calculations related to the distributions
  $\mu_1$ and $\mu_0$ take constant time, or alternatively that the
  algorithm has access to an oracle which preforms them in constant
  time. Specifically, we define below (definition~\ref{def:X}) the
  log-likelihood ratio $X = d\mu_1/d\mu_0$ and its conditional
  distributions $\nu_0$ and $\nu_1$, and assume that their cumulative
  distribution functions can be calculated in constant time. Then we
  show that the agents' computations are tractable:
  \begin{theorem}
    \label{thm:tractability}
    Fix $\mu_0$ and $\mu_1$, and let $n$ be the number of
    agents. Assume that $X$, as well as the cumulative distribution
    functions of $\nu_0$ and $\nu_1$, can be calculated in constant
    time.  Then there exists an algorithm with running time $O(nt)$,
    which, given $i$'s private signal $\psignal_i$ and the votes
    $\votevec^{t-1}=\{\vote_j(t'):\: j \in [n],\, t'<t\}$, calculates
    $\vote_i(t)$, agent $i$'s vote at time $t$.
  \end{theorem}

  We in fact provide a simple algorithm that performs this
  calculation.
\end{itemize}

\subsection{Comparison to majority voting}
Apart from being computationally easier, majority voting is inferior
in every one of the above senses. In particular, it doesn't aggregate
information as well as repeated voting until consensus, and may not
have the {\bf asymptotic learning} property. Consider the following
example: A committee has to decide whether or not to accept a
candidate for a position. Each member of the committee interviews the
candidate and forms an opinion. Now, assume that a good candidate will
make a favorable impression nine times of out ten, whereas a bad
candidate will make a favorable impression six times out of ten (being
good enough to have made it to the interview stage).  In this case,
with overwhelming probability (i.e., with probability that tends to
one as the size of the committee increases), when the candidate is
bad, about sixty percent of the committee members will still have a
good impression, and consequently a decision by majority will come to
the wrong decision, namely that the candidate is good.

This flaw is rectified by a second vote: after seeing the results of
the first round of voting, the committee members will realize that too
few of them had a good impression, and will vote against the bad
candidate in the second round. Indeed, we prove below that asymptotic
learning is always achieved after two voting rounds. This suggests
that in situations where voting until convergence to consensus is
impractical, it may be still be beneficial to have more than one round
of voting. Note that there exist other mechanisms that assure
efficient aggregation of information. For example, Gerardi and
Yariv~\cite{gerardi2007deliberative} show that adding a ``cheap talk''
deliberation phase before a strategic majority vote can also lead to
efficient aggregation\footnote{Gerardi and Yariv consider strategic
  agents and show that adding a deliberation phase can lead to equilibria
  in which information is efficiently aggregated. Essentially, the
  agents reveal their private signals and then all vote
  identically. Their work, by its nature, does not consider
  computational issues.}.

Another characteristic advantage of consensus voting is that the
strengths of the participants' convictions counts. Consider a
situation in which each agent's private signal is, with high
probability, independent of the state of the world, but with some
probability provides very convincing evidence. While a single
agent possessing the said ``smoking gun'' would have little
impact in a majority vote, his or her insistence in subsequent rounds
would convey the weight of the evidence to the rest of the group.

\subsection{Asymptotic learning vs.\ optimal aggregation of
  information}
A stronger notion than asymptotic learning is that of ``optimal
aggregation of information''. This would describe the case that the
vote that the agents eventually converge to is equal to the vote that
would be cast by a social planner who has access to all the agents'
private signals, i.e., $\argmax_{s \in
  \{0,1\}}\CondP{S=s}{\psignal_1,\ldots,\psignal_n}$.  We state here
without proof that this is not the case in our model.

This question is related to that of monotonicity. Let $w_1,\ldots,w_n$
be such that when $W_i=w_i$ then the agents all converge to 1, and let
$w'_1,\ldots,w'_n$ be such that $\CondP{S=1}{\psignal_i=w'_i} \geq
\CondP{S=1}{\psignal_i=w_i}$. Is it necessarily the case that setting
$W_i=w'_i$ would also result in agreement on 1? Perhaps surprisingly,
it is possible to construct examples in which this is not the case. A
consequence is that this model does not display optimal aggregation of
information, since $\argmax_{s \in
  \{0,1\}}\CondP{S=s}{\psignal_1,\ldots,\psignal_n}$ is monotonic in
the sense described above.

\subsection{A note on uniform priors}
We assume that the agents' prior is uniform, i.e.,
$\P{S=1}=\P{S=0}=\half$. We make this choice to simplify our notation
and make the article easier to read, although our results can be
easily extended to biased priors. For this extension an additional
requirement is needed: it is not enough that $\mu_0 \neq \mu_1$, since
it may be the case that for any value of $\psignal_i$ it holds that
$\CondP{S=1}{\psignal_i} > \half$. For example, let the prior be such
that $\P{S=1}=0.9$, and let each private signals $\psignal_i$ equal
$S$ with probability $0.51$ and equal $1-S$ with probability
$0.49$. Then for any value of $\psignal_i$ it will be the case that
$\vote_i(1)=1$. In this case, although consensus will be reached
immediately, there will be no asymptotic learning.

Hence for general priors the requirement is that $\mu_0,\mu_1$ be such
that $\P{\vote_i(1)=1} > 0$ and $\P{\vote_i(1)=0} > 0$. Given this,
the results above and the ideas of the proof below apply equally.

\section{Proofs}

Before proving our theorems we make some additional definitions. We
start by defining the log-likelihood ratio $x$ and its conditional
distributions $\nu_0$ and $\nu_1$.
\begin{definition}
  \label{def:X}
  Let $x:\R \to \R$ be given by
  \begin{align}
    \label{eq:x-def}
    x(\omega) = \log \frac{d\mu_1}{d\mu_0}(\omega).
  \end{align}

  Let $\nu_0$ be the distribution of $x(A)$ when $A \sim \mu_0$ and
  let $\nu_1$ be the distribution of $x(A)$ when $A \sim \mu_1$.
\end{definition}

Note that if $(S,\psignal) \sim \half \delta_0\otimes \mu_0 + \half
\delta_1\otimes \mu_1$, and $X=x(\psignal)$ then
\begin{align*}
  x(\psignal) = \log
  \frac{\CondP{\psignal}{S=1}}{\CondP{\psignal}{S=0}},
\end{align*}
i.e., $x(\psignal)$ is the {\em log-likelihood ratio} of $S$ given
$\psignal$.

In the proofs that follow we denote 
\begin{align*}
  X_i = x(\psignal_i)
\end{align*}
agent $i$'s private {\em log-likelihood ratio}. The advantage of
log-likelihood ratios is that they are additive for conditionally
independent signals.

In our analysis below an event that we often encounter is $a < X_i
\leq b$, and hence the following definition will be useful.
\begin{definition}
  \label{def:xab}
  Let $x : \R^2 \to \R$ be given by
  \begin{align}
    \label{eq:x-ab}
    x(a,b) = \log \frac{\mu_1(a < \omega \leq b)}{\mu_0(a < \omega \leq b)}.
  \end{align}
\end{definition}

Note that if $(S,\psignal) \sim \half \delta_0\otimes \mu_0 + \half
\delta_1\otimes \mu_1$, and  $X=x(\psignal)$ then
\begin{align*}
  x(a,b) = \log \frac{\CondP{a < X \leq b}{S=1}}{\CondP{a < X \leq
      b}{S=0}},
\end{align*}
i.e., $x(a,b)$ is the log-likelihood ratio of $S$ given that $X$ is
between $a$ and $b$.

The following claim follows by application of Bayes' law to
Eq.~\ref{eq:x-ab}. It follows from this claim that if the cumulative
distribution functions of $\nu_0$ and $\nu_1$ can be calculated in
constant time then so can $x(\cdot,\cdot)$. In what follows we will
need to use the fact that $x(\cdot,\cdot)$ can be efficiently
calculated.
\begin{claim}
  \begin{align}
    \label{eq:x-ab-calc}
    x(a,b) = \log \frac{\nu_1(X \leq b)-\nu_1(X \leq a)}
    {\nu_0(X \leq b)-\nu_0(X \leq a)}.
  \end{align}
\end{claim}

We shall also need the following easy claim in some of the proofs
below. It can be stated informally as ``the log-likelihood ratio of
the log-likelihood ratio is the log-likelihood ratio''.
\begin{claim}
  \label{claim:X-llr-X}
  Let $\mu_0,\mu_1$ be such that $d\mu_1/d\mu_0$ exists and is non-zero
  for all $\omega$, and let
  \begin{align*}
    x(\omega) = \log \frac{d\mu_1}{d\mu_0}(\omega).
  \end{align*}
  Let $\nu_0$ be the distribution of $x=x(\psignal)$ when $\psignal
  \sim \mu_0$, and let $\nu_1$ be the distribution of $x=x(\psignal)$
  when $\psignal \sim \mu_1$. Then $d\nu_1/d\nu_0$ exists and
  \begin{align}
    \label{eq:X-llr-X}
    x = \log \frac{ d \nu_1}{d \nu_0} (x).    
  \end{align}
\end{claim}
\begin{proof}
  Since $\mu_1$ and $\mu_0$ are absolutely continuous with respect to
  each other, it follows from the fact that $x(\psignal)$ is a
  function of $\psignal$ that $\nu_1$ and $\nu_0$ are also absolutely
  continuous with respect to each other, and so $d\nu_1/d\nu_0$ exists
  and is non-zero.

  Let $(S,\psignal)\sim \half
  \delta_0\otimes\mu_0+\half\delta_1\otimes\mu_1$, and denote
  $X=x(A)$. Let
  \begin{align*}
    M = \CondP{S=1}{\psignal}=\CondE{S}{\psignal}.
  \end{align*}
  Then $X = \log (M/(1-M))$.  By the law of total expectation we have
  that
  \begin{align*}
    \CondE{S}{M} = \CondE{\CondE{S}{M,\psignal}}{M}.
  \end{align*}
  Since $M$ is a function of $\psignal$ then $\CondE{S}{M,\psignal} =
  \CondE{S}{\psignal} = M$ and it follows that
  \begin{align*}
    \CondE{S}{M} = \CondE{M}{M} = M.
  \end{align*}
  Now from the definition of $X$ it follows that $X=\log(M/(1-M))$ and
  therefore there is a one-to-one correspondence between $X$ and
  $M$. Hence $\CondE{S}{X} = \CondE{S}{M} = M$. Therefore
  \begin{align*}
    \log\frac{\CondP{S=1}{X}}{\CondP{S=0}{X}} = \log\frac{M}{1-M}=X.
  \end{align*}
  But the l.h.s.\ of this equation is by Bayes' law equal to
  \begin{align*}
    \log\frac{\CondP{X}{S=1}}{\CondP{X}{S=0}},
  \end{align*}
  which is equal to $\frac{d\nu_1}{d\nu_0}(X)$, and thus we have that
  \begin{align*}
    X = \log\frac{d\nu_1}{d\nu_0}(X).
  \end{align*}
\end{proof}

\subsection{Tractability}

The key observation behind our tractability proof is the
following. Let $i$ and $j$ be two agents. Since $j$ knows all that $i$
knows except $\psignal_i$, then even before $i$ votes, agent $j$ can
know for which values of $\psignal_i$ agent $i$ would vote 1, and for
which it would vote 0. In fact, we show below that there is a bound
(that $j$ can calculate) such that if $X_i$ is below that bound then
$i$ will vote 0 and otherwise $i$ will vote 1.

Thus at each voting round, $j$ gains either a lower bound or an upper
bound on $X_i$. What we in fact show is that calculating a lower bound
$A_i(t)$ and an upper bound $B_i(t)$, on {\em all} other agents'
private likelihood ratios $X_i$, is almost all that $j$ needs to do to
calculate its votes.

The following is the definition of these lower bounds $A_i(t)$ and
upper bounds $B_i(t)$. 

While Eqs.~\ref{eq:A-def} and~\ref{eq:B-def} may seem mysterious, in
what follows we prove that these definitions indeed correspond to the
intuition provided above. To somewhat elucidate Eq.~\ref{eq:A-def}
(and similarly for Eq.~\ref{eq:B-def}) we note that that every time
agent $i$ votes 1 the other agents learn a lower bound on $X_i$; the
best lower bound is the maximum of all these, and therefore we take
the maximum over all the times that $i$ voted 1. Each of these bounds
depends on the bounds that $i$ learned on the other $X_j$'s in the
previous bounds, and therefore the terms $a_j(\bar{v}^{t'-1})$ and
$b_j(\bar{v}^{t'-1})$ appear there.

\begin{definition}
  \label{def:ab}
  Let $v_i(t) \in \{0,1\}$ for $i \in [n]$ and $t \in \N$. Similarly
  to definition~\ref{def:voting}, let $\bar{v}^t=\{v_i(t'):\: i \in
  [n],\, t' \leq t\}$ denote an element of $\{0,1\}^{nt}$.  For $i \in
  [n]$ and $t \geq 0$, let $a_i : \{0,1\}^{nt} \to \R$ and $b_i :
  \{0,1\}^{nt} \to \R$ be the functions recursively defined by
  \begin{align}
    \label{eq:A-def}
    a_i(\bar{v}^t) = \max_{t' \leq t \mbox{ s.t. } v_i(t')=1}\left\{-\sum_{j \neq i} x\left(a_j(\bar{v}^{t'-1}),\,b_j(\bar{v}^{t'-1})\right)\right\}.
  \end{align}
  and
  \begin{align}
    \label{eq:B-def}
    b_i(\bar{v}^t) = \min_{t' \leq t \mbox{ s.t. } v_i(t')=0}\left\{-\sum_{j \neq i} x\left(a_j(\bar{v}^{t'-1}),\,b_j(\bar{v}^{t'-1})\right)\right\}.
  \end{align}
  where the minimum (resp., maximum) over the empty set is taken to be
  infinity (resp., minus infinity).

  Let $A_i(t)$ and $B_i(t)$ be the random variables defined by
  \begin{align*}
    A_i(t) = a_i(\votevec^t)
  \end{align*}
  and
  \begin{align*}
    B_i(t) = b_i(\votevec^t).
  \end{align*}
\end{definition}
Note that $A_i(t)$ is non-decreasing and $B_i(t)$ is non-increasing,
in $t$. 

Since, as we show below, $A_i(t)$ and $B_i(t)$ are lower and upper
bounds on $X_i$ at time $t$, we shall need to often refer to
$x(A_i(t),B_i(t))$, and hence denote
\begin{align}
  \label{eq:x-hat-def}
  \hat{X}_i(t) = x\big(A_i(t),B_i(t)\big) = x\big(a_i(\votevec^t),b_i(\votevec^t)\big).
\end{align}

Recall that $\vote_i(t)$, agent $i$'s vote at time $t$, depends on whether
or not $\CondP{S=1}{\psignal_i,\votevec^{t-1}}$ is greater than half or
not:
\begin{align*}
  \vote_i(t) =
  \begin{cases}
    1 & \mbox{if } \CondP{S=1}{\psignal_i,\votevec^{t-1}} > 1/2 \\
    0 & \mbox{otherwise} \\
  \end{cases}
\end{align*}
Let
\begin{align}
  \label{eq:Y-def}
  \hat{Y}_i(t) = \log\frac{\CondP{S=1}{\psignal_i,\votevec^{t-1}}}{\CondP{S=0}{\psignal_i,\votevec^{t-1}}}.
\end{align}
Then
\begin{align}
  \label{eq:V-eq-Y}
  \vote_i(t)=1 \; \mbox{ iff } \; \hat{Y}_i(t) > 0.
\end{align}

\begin{theorem}
  For all $i \in [n]$ and $t>0$ it holds that
  \begin{align}
    \label{eq:boundary-sufficient}
    \hat{Y}_i(t) = X_i+\sum_{j \neq i}\hat{X}_j(t-1),
  \end{align}
\end{theorem}
\begin{proof}
  We prove by induction on $t$. The basis $t=1$ follows simply from
  the definitions; since $\votevec^0$ is empty (the agents only start
  voting at $t=1$) then $A_i(0)=-\infty$ and $B_i(0)=\infty$ for all
  $i \in [n]$, and so
  \begin{align*}
    \hat{X}_i(0) = x(A_i(0),B_i(0)) = x(-\infty,\infty) = 0,
  \end{align*}
  by the definition of $x(\cdot,\cdot)$ (Eq.~\eqref{eq:x-ab}). Another
  consequence of the fact that $\votevec^0$ is empty is that
  \begin{align*}
    \hat{Y}_i(1) = \log\frac{\CondP{S=1}{\psignal_i}}{\CondP{S=0}{\psignal_i}}
  \end{align*}
  and so for $t=1$ the statement of the theorem
  (Eq.~\eqref{eq:boundary-sufficient}) reduces to
  \begin{align*}
    \log\frac{\CondP{S=1}{\psignal_i}}{\CondP{S=0}{\psignal_i}} = X_i,
  \end{align*}
  which is precisely the definition of $X_i=x(\psignal_i)$.

  Assume the statement holds for all $t' < t$ and all $i \in
  [n]$. We will show that it holds for $t$ and all $i$.  Since, as we
  note above, $\vote_i(t')=1$ iff $\hat{Y}_i(t') > 0$, then by the
  inductive assumption we have that
  \begin{align}
    \label{eq:v-ind}
    \vote_i(t') = 1 \; \mbox { iff } \; X_i+\sum_{j \neq
      i}\hat{X}_j(t'-1) > 0
  \end{align}
  or
  \begin{align*}
    \vote_i(t') &= \ind{-\sum_{j \neq i}\hat{X}_j(t'-1)
      < X_i}\\
    &= \ind{-\sum_{j \neq
        i}x\big(a_j(\votevec^{t'-1}),b_j(\votevec^{t'-1})\big) <
      X_i},
  \end{align*}
  where the second equality follows by substituting the definition of
  $\hat{X}_i(t')$. Hence $\vote_i(t')$ is equivalent to either a lower
  bound (if it equal to 1) or upper bound (if it is equal to 0) on
  $X_i$.

  Therefore the event $\votevec^{t'}=\bar{v}^{t'}$ is equal to the
  event that
  \begin{align*}
    X_i > -\sum_{j \neq
      i}x\big(a_j(\bar{v}^{t'-1}),b_j(\bar{v}^{t'-1})\big)
  \end{align*}
  for all $i$ and $t'$ such that $v_i(t') = 1$ and
  \begin{align*}
    X_i \leq -\sum_{j \neq
      i}x\big(a_j(\bar{v}^{t'-1}),b_j(\bar{v}^{t'-1})\big)
  \end{align*}
  for all $i$ and $t'$ such that $v_i(t') = 0$. Equivalently, for all
  $i \in [n]$:
  \begin{align*}
    \max_{t' \leq t,\,
      v_i(t')=1}\left\{-\sum_{j \neq i}x\big(a_j(\bar{v}^{t'-1}),b_j(\bar{v}^{t'-1})\big)\right\}
    < X_i \leq
    \min_{t' \leq t, \,
      v_i(t')=0}\left\{-\sum_{j \neq i}x\big(a_j(\bar{v}^{t'-1}),b_j(\bar{v}^{t'-1})\big)\right\}
  \end{align*}
  Substituting the definitions of $a_i$
  (Eq.~\eqref{eq:A-def}) and $b_i$ (Eq.~\eqref{eq:B-def}), this event is
  equal to the event
  \begin{align}
    \label{eq:X-in-ab}
    a_i(\bar{v}^{t'}) < X_i \leq b_i(\bar{v}^{t'}),
  \end{align}
  for all $i \in [n]$ (note that this means that $A_i(t') < X_i \leq
  B_i(t')$ for all $i$ and $t'$). Therefore
  \begin{align*}
    \CondP{S=s}{\votevec^{t'}=v^{t'}} = \CondP{S=s}{a_i(\bar{v}^{t'}) < X_i \leq
      b_i(\bar{v}^{t'}) \mbox{ for } i \in [n]},
  \end{align*}
  and also
  \begin{align*}
    \CondP{S=s}{\psignal_i=\omega,\votevec^{t'}=\bar{v}^{t'}} = \CondP{S=s}{\psignal_i=\omega,
      a_j(v^{t'}) < X_j \leq b_j(\bar{v}^{t'}) \mbox{ for } i \neq j}.
  \end{align*}
  Hence
  \begin{align*}
    \lefteqn{\log\frac{\CondP{S=1}{\psignal_i=\omega,\votevec^{t-1}=\bar{v}^{t-1}}}{\CondP{S=0}{\psignal_i=\omega,\votevec^{t-1}=\bar{v}^{t-1}}}}\\
    &= \log\frac {\CondP{S=1}{\psignal_i=\omega, a_j(\bar{v}^{t-1})
        < X_j \leq b_j(\bar{v}^{t-1}) \mbox{ for } i \neq j}}
     {\CondP{S=0}{\psignal_i=\omega, a_j(\bar{v}^{t-1}) <
        X_j \leq b_j(\bar{v}^{t-1}) \mbox{ for } i \neq j}}.
  \end{align*}

  Again invoking Bayes' law, we have that
  \begin{align*}
    \lefteqn{\log\frac{\CondP{S=1}{\psignal_i=\omega,\votevec^{t-1}=\bar{v}^{t-1}}}{\CondP{S=0}{\psignal_i=\omega,\votevec^{t-1}=\bar{v}^{t-1}}}}\\
    &= \log \left(
      \frac {\CondP{\psignal_i=\omega}{S=1}}
      {\CondP{\psignal_i=\omega}{S=0}} \prod_{j \neq i}\frac
      {\CondP{a_j(\bar{v}^{t-1}) < X_j \leq b_j(\bar{v}^{t-1})}{S=1}}
      {\CondP{a_j(\bar{v}^{t-1}) < X_j \leq b_j(\bar{v}^{t-1})}{S=0}}
      \right)
  \end{align*}
  since the private signals are independent, conditioned on
  $S$. Substituting the definition of $x(\omega)$ (Eq.~\eqref{eq:x-def})
  and the definition of $x(\cdot, \cdot)$ (Eq.~\eqref{eq:x-ab}) yields
  \begin{align}
    \label{eq:llhr-decomp}
    \log\frac{\CondP{S=1}{\psignal_i=\omega,\votevec^{t-1}=\bar{v}^{t-1}}}{\CondP{S=0}{\psignal_i=\omega,\votevec^{t-1}=\bar{v}^{t-1}}}
    &= x(\omega)+\sum_{j \neq i}x\big(a_j(\bar{v}^{t-1}),
    b_j(\bar{v}^{t-1})\big).
  \end{align}
  Finally, since
  \begin{align*}
    \hat{Y}_i(t) =
    \log\frac{\CondP{S=1}{\psignal_i,\votevec^{t-1}}}{\CondP{S=0}{\psignal_i,\votevec^{t-1}}},
  \end{align*}
  then by Eq.~\eqref{eq:llhr-decomp}
  \begin{align*}
    \hat{Y}_i(t) =
    x(\psignal_i)+\sum_{j \neq i}x\big(a_j(\votevec^{t-1}),
    b_j(\votevec^{t-1})\big),
  \end{align*}
  and the theorem follows by substituting $X_i=x(\psignal_i)$ and $\hat{X}_j(t-1)=x(a_j(\votevec^{t-1}),b_j(\votevec^{t-1}))$.

\end{proof}

We are now ready to prove our main theorem for this subsection.
\begin{theorem}[Thm.~\ref{thm:tractability}]
  Fix $\mu_0$ and $\mu_1$, and let $n$ be the number of agents. Assume
  that $X$, as well as the cumulative distribution functions of
  $\nu_0$ and $\nu_1$, can be calculated in constant time.  Then there
  exists an algorithm with running time $O(nt)$, which, given $i$'s
  private signal $\psignal_i$ and the votes
  $\votevec^{t-1}=\{\vote_j(t'):\: j \in [n],\, t'<t\}$, calculates
  $\vote_i(t)$, agent $i$'s vote at time $t$.
\end{theorem}

\begin{proof}
  By Eq.~\eqref{eq:V-eq-Y} we have that $\vote_i(t)$ is a simple
  function of $\hat{Y}_i(t)$. By Theorem~\ref{eq:boundary-sufficient}
  above, $\hat{Y}_i(t)$ can be calculated in $O(n)$ by adding $X_i$
  (which we assume can be calculated in constant time given
  $\psignal_i$) to the sum over $j \neq i$ of $x(A_j(t), B_j(t))$. By
  Eq.~\eqref{eq:x-ab-calc}, $x(a, b)$ can be calculated in constant
  time, assuming the cumulative distribution functions of $\nu_0$ and
  $\nu_1$ can be calculated in constant time.

  We have therefore reduced the problem to that of calculating
  $A_j(t)=a_j(\votevec^t)$ and $B_j(t)=b_j(\votevec^t)$. However, the
  definitions of $a_j$ and $b_j$ (Eqs.~\eqref{eq:A-def}
  and~\eqref{eq:B-def}) are in fact simple recursive rules for
  calculating $a_j(\bar{v}^t)$ and $b_j(\bar{v}^t)$ for all $j \in
  [n]$, given $a_j^{t-1}(\bar{v}^{t-1})$ and
  $b_j^{t-1}(\bar{v}^{t-1})$ for all $j \in [n]$: it follows directly
  from Eqs.~\eqref{eq:A-def} and~\eqref{eq:B-def} that
  \begin{align*}
    a_i(\bar{v}^t) =
    \begin{cases}\max\left\{a_i(\bar{v}^{t-1}),
        \sum_{j \neq i}
        x\left(a_j(\bar{v}^{t-1}),\,b_j(\bar{v}^{t-1})\right)\right\}
      & \mbox{ if } v_i(t)=1\\
      a_i(\bar{v}^{t-1}) & \mbox{ otherwise}
      \end{cases},
  \end{align*}
  with an analogous equation for $b_i(\bar{v}^t)$.

  Note that the sum $\sum_{j \neq i}
  x\left(a_j(\bar{v}^{t-1}),\,b_j(\bar{v}^{t-1})\right)$ needn't be
  calculated from scratch for each $i$; one can rather sum over all $j
  \in [n]$ once and subtract the appropriate term for each $i$.  Hence
  calculating $a_j(\bar{v}^t)$ and $b_j(\bar{v}^t)$ (for all $j$)
  given their predecessors takes $O(n)$, and the entire recursive
  calculation takes $O(nt)$.
\end{proof}

\subsection{Unanimity}
Our model is a special case of that of Gale and
Kariv~\cite{GaleKariv:03}. They show a weak agreement result: namely,
that if the votes of two agents converge, and if the agents are not
indifferent at the limit $t \to \infty$, then they converge to the
same vote. We prove the strongest possible agreement result: consensus
is reached with probability one, i.e., the agents almost always all
converge to the same vote.

Before proving the theorem we prove some standard claims. Recall the
definition of $x(\cdot,\cdot)$ (Eq.~\eqref{eq:x-ab}):
\begin{align*}
  x(a,b) = \log \frac{\CondP{S=1}{a < X \leq b}}{\CondP{S=0}{a
      < X \leq b}}.
\end{align*}
\begin{claim}
  \label{claim:xab-log-expect-exp}
  Let $a,b$ be such that $x(a,b)$ is well defined (i.e., $\CondP{a < X
    \leq b}{S=0} > 0$). Then
  \begin{align}
    \label{eq:x-exp-ex}
    x(a,b) = \log \CondE{e^X}{a < X \leq b, S=0}.
  \end{align}
\end{claim}
\begin{proof}
  By Bayes' law we have that
  \begin{align*}
    x(a,b)= \log \frac{\CondP{a < X \leq b}{S=1}}{\CondP{a < X \leq
        b}{S=0}},
  \end{align*}
  Substituting the conditional distributions of $X$ yields
  \begin{align*}
    x(a,b) = \log \frac{ \int_a^bd\nu_1(x)}{\int_a^bd\nu_0(x) }.
  \end{align*}
  By Claim~\ref{claim:X-llr-X}
  \begin{align}
    \label{eq:x-is-x}
    x=\log\frac{d\nu_1}{d\nu_0}(x),
  \end{align}
  and so we have that
  \begin{align*}
    x(a,b) = \log \frac{\int_a^b\frac{d\nu_1}{d\nu_0}(x)d\nu_0(x)}
    {\int_a^bd\nu_0(x) } = \log \frac{\int_a^be^xd\nu_0(x)}
    {\int_a^bd\nu_0(x) }.
  \end{align*}
  Recalling that $\nu_0$ is the distribution of $X$ conditioned on $S=0$ we
  have that
  \begin{align*}
    x(a,b) = \log \CondE{e^X}{a < X \leq b, S=0}.
  \end{align*}
\end{proof}
Recall that we assume that the distribution of $X$ is non-atomic
(definition~\ref{def:assumptions}). Hence the following claim is a
consequence of Eq.~\eqref{eq:x-exp-ex} above, by a standard argument
that we omit.
\begin{claim}
  \label{claim:xab-increasing}
  $x(a,b)$ is non-decreasing and continuous in $a$ and in $b$.
\end{claim}
The following claims follows directly from Eq.~\eqref{eq:x-exp-ex} above.
\begin{claim}
  \label{claim:xab-in-ab}
  Let $a,b$ be such that $x(a,b)$ is well defined (i.e., $\CondP{a < X
    \leq b}{S=0} > 0$). Then $a < x(a,b) < b$, assuming the
  distribution of $X$ is non-atomic.
\end{claim}
\begin{proof}
  By Eq.~\eqref{eq:x-exp-ex} we have that
  \begin{align*}
  e^{x(a,b)} = \CondE{e^X}{e^a < e^X \leq e^b, S=0},
  \end{align*}
  and so $e^a < e^{x(a,b)} \leq e^b$. Since we assume the distribution
  of $X$ is non-atomic (definition~\ref{def:assumptions}) then
  \begin{align*}
    \CondE{e^X}{e^a < e^X \leq e^b, S=0} < e^b,
  \end{align*}
  and so $e^a < e^{x(a,b)} < e^b$ and the claim follows.
\end{proof}

We now show a condition for unanimity. We will later prove that
unanimity occurs w.p.\ 1 by showing that this condition
eventually applies, w.p.\ 1.
\begin{lemma}
  \label{lemma:unanimity-cond}
  If
  \begin{align}
    \label{eq:unanimity-cond}
    \sum_i \left(B_i(t)-A_i(t)\right) < \left|\sum_i X_i\right|
  \end{align}
  then there exists a $\vote$ such that $\vote_i(t') = \vote$ for all
  $i$ and $t' > t$. I.e., unanimity is reached at time $t$.
\end{lemma}
\begin{proof}
  We first note that, since $A_i(t)$ is non-decreasing and $B_i(t)$ is
  non-increasing then if Eq.~\ref{eq:unanimity-cond} holds at time $t$
  then it also holds at all times $t' > t$.

  Now, recall that $\hat{X}_i(t) = x(A_i(t),B_i(t))$. By
  Claim~\ref{claim:xab-in-ab} we have that $A_i(t) < \hat{X}_i(t) \leq
  B_i(t)$. From Eq.~\eqref{eq:X-in-ab} it follows that the same holds
  for $X_i$ too: $A_i(t) < X_i \leq B_i(t)$. Hence $|\hat{X}_i(t)-X_i|
  \leq B_i(t) - A_i(t)$ and for all $i \in [n]$ we have that 
  \begin{align}
    \label{eq:sum-less-sum}
    \left|\sum_{j \neq i}\left(\hat{X}_j(t)-X_j\right)\right| \leq
    \sum_j\left(B_j(t) - A_j(t)\right).
  \end{align}

  Recall that
  \begin{align*}
    \hat{Y}_i(t+1) = X_i+\sum_{j \neq i}\hat{X}_j(t),
  \end{align*}
  and so
  \begin{align*}
    \hat{Y}_i(t+1) -\sum_jX_j= \sum_{j \neq i}\left(\hat{X}_j(t)-X_j\right).
  \end{align*}
  Therefore by Eq.~\eqref{eq:sum-less-sum} we have that 
  \begin{align*}
    \left|\hat{Y}_i(t+1) -\sum_jX_j\right| \leq \sum_j\left(B_j(t)
      - A_j(t)\right).
  \end{align*}
  By the theorem hypothesis this implies that
  \begin{align*}
    \left|\hat{Y}_i(t+1) -\sum_jX_j\right| \leq \left|\sum_jX_j\right|.
  \end{align*}

  Hence $\hat{Y}_i(t+1)$ and $\sum_jX_j$ have the same sign. Since
  $\vote_i(t+1)=1$ iff $\hat{Y}_i(t+1) > 0$ (Eq.~\eqref{eq:V-eq-Y}) then
  we have shown that at time $t+1$ all agents vote identically. Since
  if Eq.~\ref{eq:unanimity-cond} holds for time $t$ then it
  also holds for time $t+1$ then we've shown that for all $t' > t$ the
  agents will agree in every round. It remains to show that they don't
  all change their opinion, as a group.

  Now, if the agents all vote 1 at time $t+1$ then, by the definition
  of $A_i(t)$ and $B_i(t)$, it holds that $B_i(t+2)=B_i(t+1)$ and
  $A_i(t+2) \geq A_i(t+1)$. Since by Claim~\ref{claim:xab-increasing}
  $x(a,b)$ is non-decreasing in $a$, then we have that $\hat{X}_i(t+2)
  \geq \hat{X}_i(t+1)$ for all $i$, and so $\hat{Y}_i(t+2) \geq
  \hat{Y}_i(t+1)$ for all $i$. Hence the agents will all vote 1 at
  time $t+2$. The same argument applies when all the agents vote 0 at
  time $t+1$, and the proof follows by induction on $t$.
\end{proof}

We make another definition before proceeding to prove the main theorem
of this subsection. Recall the definitions of $a_i$, $b_i$, $A_i$ and
$B_i$:
\begin{align*}
  a_i(\bar{v}^t) = \max_{t' \leq t \mbox{ s.t. } v_i(t')=1}\left\{-\sum_{j \neq i} x\left(a_j(\bar{v}^{t'-1}),\,b_j(\bar{v}^{t'-1})\right)\right\}
\end{align*}
and
\begin{align*}
  b_i(\bar{v}^t) = \min_{t' \leq t \mbox{ s.t. } v_i(t')=0}\left\{-\sum_{j \neq i} x\left(a_j(\bar{v}^{t'-1}),\,b_j(\bar{v}^{t'-1})\right)\right\}.
\end{align*}
with $A_i(t) = a_i(V^t)$ and $B_i(t) = b_i(V^t)$. As we noted above
$A_i(t)$ is non-decreasing in $t$ and $B_i(t)$ is non-increasing in
$t$. Hence they have limits which we denote by $A_i(\infty)$ and
$B_i(\infty)$. Furthermore, if as above we denote $\hat{X}_i(t) =
x\left(a_j(\votevec^{t'}),\,b_j(\votevec^{t'})\right)$ then
\begin{align}
  \label{eq:A-lim}
  A_i(\infty) = \sup_{t' \mbox{ s.t. } \vote_i(t')=1}\left\{-\sum_{j \neq i} \hat{X}_i(t')\right\}
\end{align}
and
\begin{align}
  \label{eq:B-lim}
  B_i(\infty) = \inf_{t' \mbox{ s.t. } \vote_i(t')=0}\left\{-\sum_{j \neq i} \hat{X}_i(t')\right\}
\end{align}
Note that since $\hat{X}_i(t) = x(A_i(t),B_i(t))$, and since $x(a,b)$
is a continuous function of $a$ and $b$
(Claim~\ref{claim:xab-increasing}) then
\begin{align}
  \label{eq:lim-xhat}
  \lim_{t \to \infty}\hat{X}_i(t) = x\big(A_i(\infty),B_i(\infty)\big).
\end{align}

\begin{theorem}[Thm.~\ref{thm:unanimity}]
  With probability one there exists a time $T_u$ and a vote $\vote \in
  \{0,1\}$ such that for all $t \geq T_u$ and agents $i$ it holds that
  $\vote_i(t)=\vote$.
\end{theorem}
\begin{proof}

  Assume by way of contradiction that unanimity is never reached, and
  so by Lemma~\ref{lemma:unanimity-cond} for all $t$ it holds that
  $\sum_i B_i(t)-A_i(t) \geq |\sum_iX_i|$. Then, since $B_i(t)-A_i(t)$
  is monotonically decreasing, it holds that
  \begin{align}
    \label{eq:sum_of_bounds}
    \lim_{t\to\infty} \sum_i \left(B_i(t)-A_i(t)\right) \geq \left|\sum_i X_i\right|.
  \end{align}

  Let $Z:=\lim_{t\to\infty}\sum_j\hat{X}_j(t)$.  We consider
  separately the events that $Z=0$, $Z<0$ and $Z>0$:
  \begin{enumerate}
  \item $Z=0$

    We assume (definition~\ref{def:assumptions}) that the distribution
    of $X_i$ is non-atomic, and so $\sum_i X_i \neq 0$ with
    probability $1$. Hence by Eq.~\eqref{eq:sum_of_bounds} there must
    be some agent $i$ for which
    \begin{align*}
    \lim_{t\to\infty}
    \left(B_i(t)-A_i(t)\right)=B_i(\infty)-A_i(\infty)>0.  
    \end{align*}
    Assume w.l.o.g.\ $V_i(t) = 1$ infinitely many
    times. Hence, by Eq.~\eqref{eq:A-lim}, we have that
    \begin{align*}
      A_i(\infty) &\geq -\lim_{t \to \infty}\sum_{j \neq i}
      \hat{X}_i(t)\\
      &= \lim_{t \to \infty}\hat{X}_i(t)-\lim_{t \to \infty}\sum_j \hat{X}_j(t).
    \end{align*}
    Since we assume in this case that
    $Z=\lim_{t\to\infty}\sum_j\hat{X}_j(t)=0$ then we have that
    \begin{align*}
      A_i(\infty) \geq \lim_{t \to \infty}\hat{X}_i(t).
    \end{align*}
    But since $A_i(\infty)<B_i(\infty)$ then by
    Eq.~\eqref{eq:lim-xhat} and Claim~\ref{claim:xab-in-ab} we have
    that $A_i(\infty) < \lim_{t\to\infty}\hat{X}_i(t)$, which is
    a contradiction.

    The intuition here is that when $i$ votes 1 it is revealed that
    $X_i>\hat{X}_i(t)-\sum_j\hat{X}_j(t)$. Hence if
    $\sum_j\hat{X}_j(t)$ is very small then $A_i(t)$ approaches
    $\hat{X}_i(t)$ arbitrarily closely, which is impossible if
    $A_i(t)$ is to stay well separated from $B_i(t)$.
  \item $Z > 0$

    Since unanimity is never reached then there must be some $i$ for
    which $\hat{Y}_i(t) \leq 0$ infinitely many times.  Hence by
    Eq.~\eqref{eq:B-lim} we have that
    \begin{align*}
      B_i(\infty) &\leq -\lim_{t \to \infty}\sum_{j \neq i}
      \hat{X}_i(t)\\
      &= \lim_{t \to \infty}\hat{X}_i(t)-\lim_{t \to \infty}\sum_j \hat{X}_j(t).
    \end{align*}
    Since by assumption $\lim_{t \to \infty} \sum_j \hat{X}_j(t) > 0$
    then we have that
    \begin{align*}
      B_i(\infty) < \lim_{t \to \infty}\hat{X}_i(t).
    \end{align*}
    However by Eq.~\eqref{eq:lim-xhat} and Claim~\ref{claim:xab-in-ab}
    we have that $B_i(\infty) > \lim_{t\to\infty}\hat{X}_i(t)$, which
    is a contradiction.

    In this case the intuition is that when $i$ votes 0 even though
    $\sum_j \hat{X}_j(t)$ is positive, then $B_i$ decreases by at
    least $\sum_j \hat{X}_j(t)$, which cannot continue indefinitely
    when $\lim_{t\to\infty}\sum_j \hat{X}_j(t) > 0$.

  \item $Z < 0$

    The argument here is identical to that of the previous case.

  \end{enumerate}
\end{proof}

\subsection{Monotonicity}
Since the agents base their decisions on a growing information
base, their decisions become more and more likely to be correct. We
prove this formally below, using a standard argument.

\begin{theorem}[Thm.~\ref{thm:monotonicity}]
  For all agents $i$ and times $t>1$, it holds that
  \begin{align*}
    \P{\vote_i(t)=S} \geq \P{\vote_i(t-1)=S}.
  \end{align*}
\end{theorem}

\begin{proof}
  As noted in Eq.~\eqref{eq:vote-max-prob}, $\vote_i(t)$ is the choice
  in $\{0,1\}$ that maximizes the probability of matching the state of
  the world, given $\psignal_i$ and $\votevec^{t-1}$. Let $f$ be an
  arbitrary function of $\psignal_i$ and $\votevec^{t-1}$. Then:
  \begin{align*}
    \CondP{\vote_i(t)=S}{\psignal_i,\votevec^{t-1}} \geq
    \CondP{f\big(\psignal_i,\votevec^{t-1}\big)=S}{\psignal_i,\votevec^{t-1}}.
  \end{align*}
  Since $\vote_i(t-1)$ is also a function of $\psignal_i$ and
  $\votevec^{t-1}$ then we can substitute $\vote_i(t-1)$ for
  $f(\psignal_i,\votevec^{t-1})$ in the equation above, and the theorem
  follows.
\end{proof}
Note that $\P{\vote_i(t)=S}$ is strictly larger than
$\P{\vote_i(t-1)=S}$ whenever $\P{\vote_i(t)\neq \vote_i(t-1)}$ is
positive, i.e. when the decision may change.

\subsection{Asymptotic Learning}
We show that with high probability after observing the first round of
voting all voters know the correct state of the world, and a unanimous
and correct decision is reached at the second round of voting. Note
that by the monotonicity theorem (\ref{thm:monotonicity}), this means
that the same holds for all rounds after the second round.

Before proving this theorem we will prove the following claim.
\begin{claim}
  \label{claim:X-markov}
  \begin{align*}
    \CondP{X < -C}{S=1} < e^{-C}.
  \end{align*}
\end{claim}
\begin{proof}
  Recall that by Claim~\ref{claim:X-llr-X}
  \begin{align*}
    x = \frac{d\nu_1}{d\nu_0}(x)
  \end{align*}
  and so
  \begin{align*}
    \int_{-\infty}^\infty d\nu_0(X) = \int_{-\infty}^\infty e^{-X}d\nu_1(X).
  \end{align*}
  But $\nu_0$ is a probability measure, and so $\int_{-\infty}^\infty
  d\nu_0(X)=1$. Hence
  \begin{align*}
    \CondE{e^{-X}}{S=1} = \int_{-\infty}^\infty e^{-X}d\nu_1(X) = 1.
  \end{align*}
  Therefore by the Markov bound
  \begin{align*}
    \CondP{X < -C}{S=1} = \CondP{e^{-X} > e^C}{S=1} < e^{-C}.
  \end{align*}
\end{proof}

We are now ready to prove the main theorem of this subsection.
\begin{theorem}[Thm.~\ref{thm:learning}]
  Fix $\mu_0$ and $\mu_1$, and let $n$ be the number of agents.  Then
  there exist constants $C=C(\mu_0,\mu_1)$ and $n_0=n_0(\mu_0,\mu_1)$ such
  that
  \begin{align*}
    \P{\forall i: \vote_i(2) = S} > 1-e^{-Cn}
  \end{align*}
  for all $n > n_0$.
\end{theorem}

\begin{proof}
  We shall show that there exists a constant $C=C(\mu_0,\mu_1)$ such that
  \begin{align*}
    \CondP{\forall i: \vote_i(2) = 1}{S=1} > 1-e^{-Cn}
  \end{align*}
  for all $n$ large enough. Since the same argument can be used to
  show an analogous statement for $S=0$ then this will prove the
  theorem.

  Recall that $\vote_i(t)$ is the indicator of the event $\hat{Y}_i(t)
  > 0$, where
  \begin{align*}
    \hat{Y}_i(t) =
    \log\frac
    {\CondP{S=1}{\psignal_i,\votevec^{t-1}}}
    {\CondP{S=0}{\psignal_i,\votevec^{t-1}}}.
  \end{align*}
  Invoking Bayes' law and the conditional independence of the private
  signals we get that
  \begin{align*}
    \hat{Y}_i(2) = \log \left(
      \frac {\CondP{\psignal_i}{S=1}}{\CondP{\psignal_i}{S=0}}
      \prod_{j \neq i}\frac
      {\CondP{\vote_j(1)}{S=1}}
      {\CondP{\vote_j(1)}{S=0}}
      \right)
  \end{align*}
  Denote by $N_i$ the number of agents other than $i$ who vote 1 in
  the first round:
  \begin{align*}
    N_i = |\{j \mbox{ s.t. } \vote_j(1)=1, j \neq i\}|.
  \end{align*}
  Then since $X_i = \log\frac{\CondP{\psignal_i}{S=1}}{\CondP{\psignal_i}{S=0}}$
  then we can write
  \begin{align*}
    \hat{Y}_i(2) &= X_i+N_i\log\frac
    {\CondP{X>0}{S=1}}
    {\CondP{X>0}{S=0}}
    +(n-1-N_i)\log\frac
    {\CondP{X \leq 0}{S=1}}
    {\CondP{X \leq 0}{S=0}}.
  \end{align*}
  Denote
  \begin{align*}
    \alpha_1=\CondP{X > 0}{S=1}\quad \mbox{and} \quad
    \alpha_0=\CondP{X > 0}{S=0},
  \end{align*}
  and note that $x(0,\infty) > 0$ by Claim~\ref{claim:xab-in-ab}, and
  so since $x(0,\infty) = \log\alpha_1/\alpha_0$ then $\alpha_1 \neq
  \alpha_0$.  We can now write
  \begin{align*}
    \hat{Y}_i(2) &= X_i + N_i\log\frac
    {\alpha_1}
    {\alpha_0}
    +(n-1-N_i)\log\frac
    {1-\alpha_1}
    {1-\alpha_0}.
  \end{align*}
  Denote
  \begin{align*}
    W_i = N_i\log\frac {\alpha_1} {\alpha_0} +(n-1-N_i)\log\frac
    {1-\alpha_1} {1-\alpha_0}
  \end{align*}
  so that $\hat{Y}_i(2) = X_i + W_i$. Since
  $\CondE{N_i}{S=1}=(n-1)\alpha_1$ then the conditioning $W_i$ on
  $S=1$ we get that
  \begin{align*}
    \CondE{W_i}{S=1} =
    (n-1)\alpha_1\log\frac{\alpha_1}{\alpha_0}+(n-1)(1-\alpha_1)\log\frac{1-\alpha_1}{1-\alpha_0}.
  \end{align*}
  If we denote $(n-1)D=\CondE{W_i}{S=1}$ then $D$ is the
  Kullback-Leibler divergence~\cite{kullback1951information} of two
  Bernoulli distributions with expectations $\alpha_1 \neq
  \alpha_0$. Hence $D > 0$.

  Now, conditioned on $S$ the private signals are independent, and
  hence so are the votes at round 1, since $V_i(1)$ depends on
  $\psignal_i$ only.  Therefore $W_i$, conditioned on $S=1$, is the sum
  of $n-1$ bounded independent random variables. Therefore by the
  Hoeffding bound there exists a constant $C_1$ such that
  \begin{align}
    \label{eq:w_1-bound}
    \CondP{W_i \leq \half (n-1)D}{S=1} \leq e^{-C_1 (n-1)}.
  \end{align}

  By Claim~\ref{claim:X-markov} we have that $\CondP{X_i \leq -\half
    (n-1)D}{S=1} < e^{-\half D (n-1)}$, which, together with
  Eq.~\eqref{eq:w_1-bound} and the union bound yields
  \begin{align*}
    \CondP{\vote_i(2) = 0}{S=1} < e^{-\half D (n-1)} + e^{-C_1 (n-1)}.
  \end{align*}
  Therefore, again using the union bound it follows that
  \begin{align*}
    \CondP{\exists i: \vote_i(2) = 0}{S=1} < n(e^{-\half D (n-1)} +
    e^{-C_1 (n-1)}).
  \end{align*}
  Finally, it follows that for $n$ large enough there exists a
  constant $C$ such that
  \begin{align*}
    \CondP{\forall i: \vote_i(2) = 1}{S=1} > 1- e^{-C n}.
  \end{align*}
\end{proof}
\bibliographystyle{abbrv}      
\bibliography{all}

\end{document}